\documentclass[journal,twoside,shortpaper,9pt]{IEEEtran}

\usepackage{atbegshi,picture,ifpdf,setspace}
\usepackage[noadjust]{cite}
\usepackage[pdftex]{graphicx}
\usepackage[cmex10]{amsmath}
\usepackage[caption=false,font=footnotesize]{subfig}
\usepackage{url}
\usepackage[english]{babel}
\usepackage[utf8x]{inputenc}
\usepackage[T1]{fontenc}
\usepackage[scaled=0.91]{helvet}
\usepackage[cmintegrals,varvw,]{newtxmath}
\usepackage{bm}
\usepackage{siunitx} 
\addto\captionsenglish{}

\makeatletter
\def\ps@IEEEtitlepagestyle{%
  \def\@oddfoot{\mycopyrightnotice}%
  \def\@oddhead{\hbox{}\@IEEEheaderstyle\leftmark\hfil\thepage}\relax
  \def\@evenhead{\@IEEEheaderstyle\thepage\hfil\leftmark\hbox{}}\relax
  \def\@evenfoot{}%
}
\def\mycopyrightnotice{%
  \begin{minipage}{\textwidth}
  \centering \scriptsize
  Copyright~\copyright~2018 IEEE. Personal use of this material is permitted. Permission from IEEE must be obtained for all other uses, in any current or future media, including\\reprinting/republishing this material for advertising or promotional purposes, creating new collective works, for resale or redistribution to servers or lists, or reuse of any copyrighted component of this work in other works by sending a request to pubs-permissions@ieee.org.
  \end{minipage}
}
\makeatother

\AtBeginShipout{\AtBeginShipoutUpperLeft{%
  \put(\dimexpr\paperwidth-1cm\relax,-.55cm){\makebox[0pt][r]{
  \begin{minipage}{\textwidth}
  \centering \scriptsize 
  This is the author's version of an article that has been published in this journal. 
  Changes were made to this version by the publisher prior to publication. 
  \\[0.3ex]
  The final version of record is available at\qquad http://dx.doi.org/10.1109/TAP.2018.2800734
  \end{minipage}
  }}}%
}

\begin{document}
\title{{\fontsize{24}{26}\selectfont{Communication\rule{29.9pc}{0.675pt}}}\vspace*{0.2cm}\break\fontsize{16}{18}\selectfont
\vspace*{-0.1cm}%
A Combined Source Integral Equation\\[0.75ex]with Weak Form Combined Source Condition}
\author{
\Large Jonas Kornprobst,  \IEEEmembership{\Large Student Member, IEEE}, and Thomas F. Eibert, \IEEEmembership{\Large Member, IEEE}%
\thanks{Manuscript received December 09, 2016; revised October 05, 2017; accepted
January 22, 2018; date of this version October 17, 2020. \emph{(Corresponding author: Jonas Kornprobst.)}}%
\thanks{J. Kornprobst and T. F. Eibert are with the Chair of High-Frequency Engineering, Department of Electrical and Computer Engineering, Technical University of Munich, 80290 Munich, Germany (e-mail: j.kornprobst@tum.de; hft@ei.tum.de).}%
\thanks{Color versions of one or more of the figures in this communication are available online at http://ieeexplore.ieee.org.}%
\thanks{Digital Object Identifier 10.1109/TAP.2018.2800734}%
}
\markboth{IEEE TRANSACTIONS ON ANTENNAS AND PROPAGATION}{Kornprobst and Eibert: A Combined Source Integral Equation with Weak-Form Combined Source Condition}

\maketitle

\begin{abstract}
A combined source integral equation (CSIE) is constructed on the basis of the electric field integral equation (EFIE) to solve electromagnetic radiation and scattering problems containing perfect electrically conducting bodies. 
It is discretized with Rao-Wilton-Glisson basis functions only, for both electric and magnetic surface current densities. 
The combined source condition, which ensures the uniqueness of the solution and circumvents the interior resonance problem, is implemented as a weak form side condition. 
Compared to the common combined field integral equation, the proposed CSIE shows superior accuracy for sharp edges as well as structures with the interior resonance problem. 
Furthermore, the iterative solver convergence of the CSIE is faster than for the EFIE, which shows about the same accuracy as the CSIE. 
Results of numerical scattering simulations are presented to demonstrate the accuracy of the presented CSIE. 
\end{abstract}


%%%%%%%%%%%%%%%%%%%%%%%%%%
\section{Introduction} %%%
%%%%%%%%%%%%%%%%%%%%%%%%%%
\IEEEPARstart{S}{urface} integral equation (IE) solutions by the method of moments (MoM) are the method of choice for the solution of electromagnetic radiation and scattering problems involving perfect electrically conducting (PEC) objects. 
The fundamental surface IEs are the electric field integral equation (EFIE) and the magnetic field integral equation (MFIE), which impose boundary conditions on the tangential electric or magnetic field, respectively~\cite{peterson1998computational}. 
The EFIE is suitable for open and closed structures, whereas the MFIE works for closed structures only. 
Both the {EFIE} and the {MFIE} suffer from the interior resonance problem when applied to closed bodies, since their solutions are not unique at the corresponding interior resonance frequencies and can be superimposed by non-radiating currents~\cite{Peterson1990}. 
Combining EFIE and MFIE to the so-called combined field integral equation (CFIE) is the most popular approach to circumvent this issue~\cite{mautz1979h}. 
Since the MFIE is a second kind IE, a widely appreciated side effect of utilizing the CFIE is its improved iterative solver convergence as compared to pure first kind EFIE solutions. 

In most cases, the modeling of arbitrary surfaces is performed based on triangular facets together with the corresponding Rao-Wilton-Glisson (RWG) basis functions, which are the lowest-order divergence-conforming functions on triangular meshes~\cite{Rao_1982}. 
Formally, the RWG functions are the appropriate choice for the discretization of surface current densities. 
The testing functions should be chosen from the dual space of the range of the considered operator acting on the surface current densities~\cite{Yla-OijalaMarkkanenJarvenpaaEtAl2014}. 
However, especially with the low-order RWG functions, it is not always feasible to follow these rules and it was found that the MFIE (and also CIFE) suffers from inaccuracy problems due to the discretization and testing of the MFIE operator, especially for sharp edges and electrically small features~\cite{Gurel_2009}. 
Over the years, different types of basis and testing functions have been employed  to analyze and enhance the accuracy of the {MFIE}~\cite{Ergul_2007,Ismatullah_2009_MFIE,Peterson2008}. 
With low-order basis functions, well-conditioned accurate formulations are only possible with the curl- and div-conforming Buffa-Christiansen (BC) functions~\cite{Cools_2011}. 
However, in many situations it is still desirable to work with RWG functions only and many techniques applicable to the MFIE may not be appropriate for use with the CFIE. 
In any case, it is of course still necessary to add higher-order basis functions, if higher-order accuracy is needed. 

A very interesting interior resonance-free alternative to the CFIE is the combined source integral equation (CSIE) \cite{Bolomey_1973,Mautz1979,Rogers1985,Lee2005}, which has been known for many years, but which has never really gained a lot of popularity, most likely due to its obviously more involved numerical treatment. 
More recently, the CSIE has been revisited and discretized with different basis functions for electric and magnetic surface current densities. 
In~\cite{Yla_Oijala_2012,Ylae-OijalaKiminkiJaervenpaeae2013}, RWG and BC basis functions have been employed for electric and magnetic surface currents in scattering computations involving dielectric bodies. 
In particular, the combined source (CS) conditions have been implemented in a weak form in order to avoid the undesirable line charges at the edges of the discretization elements. 
Conceptually, the CS condition is nothing else than an impedance boundary condition (IBC), which has already been implemented in weak form with RWG functions only in~\cite{Ismat2009MoM} and later also in~\cite{YanJin2013}. 
In \cite{EibertVojvodicHansen2016}, a strong form CS condition has been utilized in the inverse source solution for antenna field transformations and a weak form CS inverse source solution is discussed in \cite{Eibert2017}. 

In this paper, the CSIE is investigated with respect to the solution of electromagnetic scattering problems involving PEC bodies. 
In order to avoid BC functions, we work with RWG basis functions only, for the discretization of both electric and magnetic surface current densities. 
The CS condition is implemented in a weak form and the additional CS equations are solved as a part of the overall equation system without explicit inversion of the Gram matrix.  
Compared to standard MFIE and CFIE solutions, this approach shows a superior accuracy for structures with sharp edges and the interior resonance problem is also well under control. 
The iterative solver convergence of the CSIE is faster than for the EFIE and both exhibit about the same accuracy, except for the  interior resonance of course, where the EFIE is known to fail.
 
In section II, the formulation and implementation of the CSIE with weak form CS condition and pure RWG discretization (basis and testing function) is presented. Thereafter in section III, its performance is demonstrated and investigated by a series of numerical results. 

%%%%%%%%%%%%%%%%%%%%%%%%%
\section{Formulation} %%%
%%%%%%%%%%%%%%%%%%%%%%%%%

\begin{figure}[tp]
 \centering
 \begin{minipage}[b]{0.49\linewidth}
  \centering
  \subfloat[\label{fig:scatterer_a}]{%
   \includegraphics[]{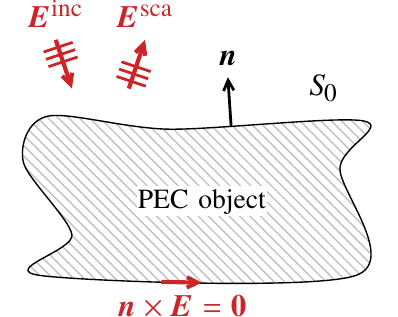}%
   \hspace*{0.25cm}%
   }%
 \end{minipage}
 \begin{minipage}[b]{0.49\linewidth}
  \raggedleft
  \subfloat[\label{fig:scatterer_b}]{%
   \hspace*{0.25cm}%
   \includegraphics[]{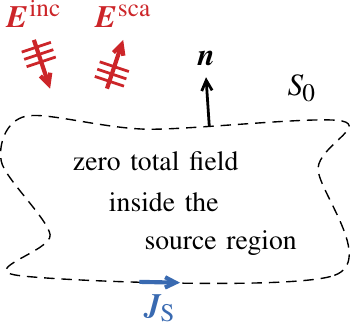}%
  }%
 \end{minipage}\\[2ex]
 \begin{minipage}[b]{0.99\linewidth}
  \centering
  \subfloat[\label{fig:scatterer_c}]{%
   \hspace*{0.25cm}%
   \includegraphics[]{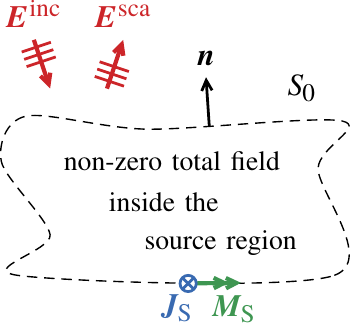}%
   \hspace*{0.25cm}%
  }%
 \end{minipage}
 \caption{PEC object homogeneously filled with the solution domain material and modeled by different choices of  equivalent sources. (a) PEC object with incident and scattered electric field. (b) PEC object, replaced by Love currents. (c) Object, replaced by CS currents.\label{fig:scatterer}} 
 \vspace*{-0.2cm}
\end{figure}

Consider the tangential component $\bm n(\bm r) \times \bm E (\bm r)$ of the electric field on the surface $S_0$ of an object~\cite{peterson1998computational}
\begin{equation}
\bm n(\bm r) \times \bm E(\bm r) = \bm n(\bm r) \times \bm E^\mathrm{inc}(\bm r) + \bm n(\bm r) \times \bm E^\mathrm{sca}(\bm r),\label{eq:scattered_field}
\end{equation}
where $\bm E^\mathrm{inc}(\bm r)$ is the incident field and $\bm E^\mathrm{sca}(\bm r)$ is the scattered field.\footnote{For radiation problems, the incident field would be directly impressed on the surface of the object, e.g.\ in form of a delta voltage source.} 
The observation point  on the surface of the object is $\bm r$ and  $\bm n$  is the unit normal vector on the surface $S_0$ pointing into the solution domain. 
This scenario, with a PEC scatterer, is depicted in Fig.~\ref{fig:scatterer}(a). 
The fields are assumed to be time harmonic with a suppressed time dependence~$\mathrm{e}^{\,\mathrm{j}\omega t}$. 
Replacing the object by equivalent electric and magnetic surface current densities $\bm J_\mathrm S$ and $\bm M_\mathrm S$, respectively, on the surface $S_0$, equation \eqref{eq:scattered_field} can be written as \cite{peterson1998computational,Mautz1979} 
\begin{IEEEeqnarray}{r}
\bm n \times \bm E(\bm r) = \bm n \times \bm E^\mathrm{inc}(\bm r) 
-\bm n \times \oiint\limits_{S_0}\bm\nabla G_0(\bm r,\bm r')\times\bm M_\mathrm S(\bm r') \mathop{}\!\mathrm{d}s'
\nonumber\\
-\frac{1}{2}\bm M_\mathrm S (\bm r) 
- \mathrm{j} k_0 Z_0 \bm n \times 
\Bigg[
\frac{1}{k_0^2}\bm\nabla\oiint\limits_{S_0}G_0(\bm r,\bm r')\bm\nabla'\cdot\bm J_\mathrm S(\bm r')\mathop{}\!\mathrm{d}s'
\nonumber\\
+\oiint\limits_{S_0}G_0(\bm r,\bm r') \bm J_\mathrm S(\bm r')\mathop{}\!\mathrm{d}s'
\Bigg],
\IEEEeqnarraynumspace\label{eq:EFIE_general}
\end{IEEEeqnarray}
 with the free space impedance $Z_0$, the free space wave number $k_0$, the source locations $\bm r'$  and the free space Green's function 
\begin{equation}
G_0(\bm r,\bm r') = \frac{\mathrm e ^{-\mathrm j k_0 \left|\bm r -\bm r'\right|}}{4\uppi\left|\bm r -\bm r'\right|}.
\end{equation} 
The common standard choice for the equivalent surface current densities is the so-called Love currents~\cite{Schelkunoff1936,Love1901}
\begin{equation}
\bm J_\mathrm S(\bm r) = \bm n \times \bm H(\bm r), \qquad \bm M_\mathrm S (\bm r)= \bm E (\bm r)\times \bm n, \label{eq:Love}
\end{equation}
which are related to the total  tangential field components on the surface of the object. 
These currents produce zero total fields inside the volume of the PEC object. 
This equivalent scenario is shown in Fig.~\ref{fig:scatterer}(b). 

Choosing the Love currents and assuming a PEC object, the condition $\bm M_\mathrm S=\bm E \times \bm n=0$ must be satisfied and converts equation \eqref{eq:EFIE_general} into the standard EFIE with electric currents only. 
In the concept of the CSIE, both current types are retained and the obvious non-uniqueness of the solution is overcome by an additional CS constraint equation. 
The PEC boundary condition $\bm n \times \bm E=0$ still has to be fulfilled and \eqref{eq:EFIE_general} is now converted into an EFIE in the form~\cite{Mautz1979}
\begin{IEEEeqnarray}{rt}
\bm n \times \bm E^\mathrm{inc}(\bm r)=\frac{1}{2}\bm M_\mathrm S(\bm r)  + \bm n \times 
\Bigg[
\oiint\limits_{S_0}\bm\nabla G_0(\bm r,\bm r')\times\bm M_\mathrm S(\bm r') \mathop{}\!\mathrm{d}s'
\nonumber\\
 +\mathrm{j} k_0 Z_0\bigg(\frac{1 }{k_0^2}\bm\nabla\oiint\limits_{S_0}G_0(\bm r,\bm r')\bm\nabla'\cdot\bm J_\mathrm S(\bm r')\mathop{}\!\mathrm{d}s'
\nonumber\\
+ \oiint\limits_{S_0}G_0(\bm r,\bm r') \bm J_\mathrm S(\bm r')\mathop{}\!\mathrm{d}s'
\bigg) \Bigg],\label{eq:EFIE}
\end{IEEEeqnarray}
which contains both types of equivalent surface current densities. 
This non-unique EFIE is augmented by an additional CS condition in the form of~\cite{Rogers1985,Mautz1979}
\begin{equation}
\bm M_\mathrm S(\bm r) = \alpha Z_\mathrm 0 \big(\bm n \times \bm J_\mathrm S(\bm r)\big),\label{eq:IBC_CS}
\end{equation}
where $\alpha$ is a weighting factor. 
With $\alpha\!=\!-1$, this condition would be nothing else than a standard Leontovich IBC for an impedance body with a surface impedance $Z_0$. 
With the common CS choice of $\alpha\!=\!1$ (or at least positive) the equivalent sources produce wave propagation into the outside solution domain and, thus, non-radiating interior cavity modes are effectively suppressed. 
In Fig.~\ref{fig:scatterer}(c), the CS solution with non-vanishing total fields in the source domain is depicted.

The choice $\alpha\!=\!0$ results in electric surface currents only, whereas $\alpha\!\rightarrow\!\infty$ leads in essence to an equation with magnetic currents only. 
This latter equation would be an IE of second kind, which contains the same operators as the MFIE. 
For the standard {CSIE}, the weighting factor is chosen to be $\alpha\!\!=\!\!1$, which  results in equal influence of electric and magnetic currents~\cite{Mautz1979}. 
Interesting to note is that the combined sources consisting of perpendicular electric and magnetic currents are nothing else than the well-known Huygens radiators, which produce a directive radiation characteristic, in our case with a null in the direction along the negative surface normal~\cite{schelkunoff1952antennas,Luk2012,Niemi2012}. 
As a result, the propagation of scattered fields through the scatterer is partly suppressed resulting in an improved conditioning of the resulting equation system, as compared to the case of the EFIE. 

In our approach, both electric and magnetic surface current densities are treated as unknowns and both are modeled by {RWG} basis functions $\bm \beta$ according to
\begin{equation}
\bm J_\mathrm S = \sum_{n=1}^N I_n \bm \beta_n, \qquad \bm M_\mathrm S = \sum_{n=1}^N V_n \bm \beta_n,
\end{equation}
where $I_n$ and $V_n$ are the $N\!+\!N$ unkown expansion coefficients. 
The electric field is tested by RWG functions, i.e.\ in its dual space. 
The {MoM} equation system resulting from testing~$\bm n \times \bm E$ in \eqref{eq:EFIE} with~$\bm n \times$RWG functions is obtained as
\begin{IEEEeqnarray}{r}
 \sum_{n=1}^N\left\{ \Bigg[-\frac{1}{2} A_{mn} +D_{mn}\Bigg] V_n +  \mathrm j Z_0 k_0\Bigg[B_{mn} +\frac{C_{mn}}{k_0^2}\Bigg]\right\} I_n=G_m^E,
 \nonumber\\
 m=1,...,N,\IEEEeqnarraynumspace
 \label{eq:EFIE_discrete}
\end{IEEEeqnarray}
with the coupling integrals $A_{mn}$, $B_{mn}$, $C_{mn}$, $D_{mn}$ and the tested incident electric field $G_m^E$~\cite{Ismat2009MoM}. 
In particular, the identity operator is evaluated as
\begin{equation}
A_{mn} =- \iint\limits_{S_m} \big(\bm n(\bm r) \times \bm \beta_n(\bm r) \big) \cdot\bm \beta_m(\bm r) \mathop{}\!\mathrm{d}s.
\end{equation}
Discretizing and testing~\eqref{eq:IBC_CS} by RWG functions results in
\begin{equation}
\sum\limits_{n=1}^N\Big\{ A_{mn}' V_n +\alpha Z_0 A_{mn} I_n\Big\}=0, \quad m=1,...,N,\label{eq:IBC_CS_weak}
\end{equation}
where the electric currents are projected on the Gram matrix~$A$ linking RWG and $\bm n \times$RWG and the magnetic currents are projected on the Gram matrix~$A'$ of RWG functions with the matrix entries defined as
\begin{equation}
A_{mn}' = \iint\limits_{S_m} \bm \beta_m(\bm r) \cdot  \bm \beta_n (\bm r)\mathop{}\!\mathrm{d}s.
\end{equation}

For the iterative solution, the discretized equation system contains the EFIE part~\eqref{eq:EFIE_discrete} and the CS side condition~\eqref{eq:IBC_CS_weak} for enforcing the unique relation between electric and magnetic currents. 

\begin{figure}[tp]
 \centering
 \begin{minipage}[b]{0.35\linewidth}
  \centering
  \subfloat[\label{fig:structures_a}]{%
   \includegraphics[]{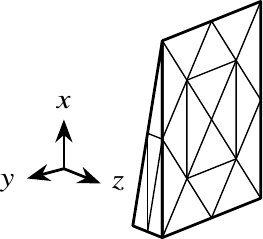}%
   }%
 \end{minipage}
 \begin{minipage}[b]{0.63\linewidth}
  \raggedleft
  \subfloat[\label{fig:structures_b}]{%
   \hspace*{0.5cm}%
   \includegraphics[]{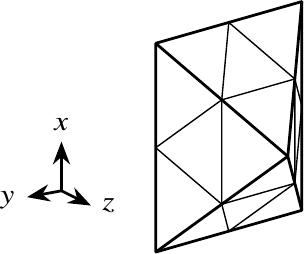}%
   \hspace*{0.5cm}%
  }%
 \end{minipage}
 \caption{Models for showing the MFIE accuracy problems, meshed with triangles for RWG basis functions. (s) Meshed wedge. (b) Meshed model of the pyramid.\label{fig:structures}} 
\end{figure}

\section{Numerical Results}

The goal of this section is to show that the new CSIE formulation with weak form CS condition and pure RWG discretization and testing leads to improved accuracy as compared to the traditional CFIE based on RWG functions. 
Since the known accuracy problems are mostly due to the MFIE,  the error due to the MFIE within the CFIE is demonstrated in particular. 
Corresponding to this is the error due to the magnetic currents within the CSIE. 
Since a pure magnetic current CSIE is not well tested, $\alpha$ in \eqref{eq:IBC_CS} is not chosen larger than 9, which is somehow similar to a CFIE with a combination parameter 0.1 (i.e.\ 10\si{\percent} EFIE and 90\si{\percent} MFIE).

Since the accuracy problem of the {MFIE} arises from objects with sharp edges and  electrically small features, the two test objects depicted in Fig.~\ref{fig:structures} are considered. 
As a larger and more realistic problem, scattering at the stealth plane Flamme is simulated. 
For more specific investigations on the good behavior at interior resonances, see the results in~\cite{KornprobstAPS}.

\subsection{Small Wedge}

A first example is a small wedge~\cite{Ismatullah_2009_MFIE}, whose structure including the mesh is shown in Fig.~\ref{fig:structures}(a). 
The height and width are $0.2\lambda$, respectively, the depth in \textit{z} direction is $0.06\lambda$. 
The discretization elements have a mean edge length of approximately $0.09\lambda$. 
The number of electric (and magnetic) current unknowns is 60, on a mesh with 40 triangles. 

For this structure, the proposed CSIE has been applied in two different ways. 
First, the CFIE with a combination parameter of $0.1$ ($0.1$-CFIE) is compared to a CSIE with approximately only magnetic currents, with $\alpha=9$ (M-CSIE). 
Secondly, the CFIE with equal weighting of EFIE and MFIE is compared to the CSIE with $\alpha\!=\!1$. 

\begin{figure}[tp]
 \begin{minipage}[b]{\linewidth}
  \centering
  \subfloat[\label{fig:RCS_wedge_a}]{%
\includegraphics{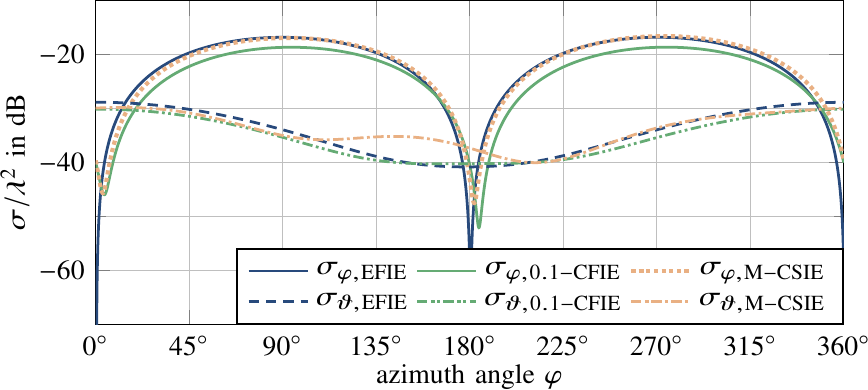}%
   }%
 \end{minipage}
 \\[2ex]
 \begin{minipage}[b]{\linewidth}
  \raggedleft
  \subfloat[\label{fig:RCS_wedge_b}]{%
\includegraphics{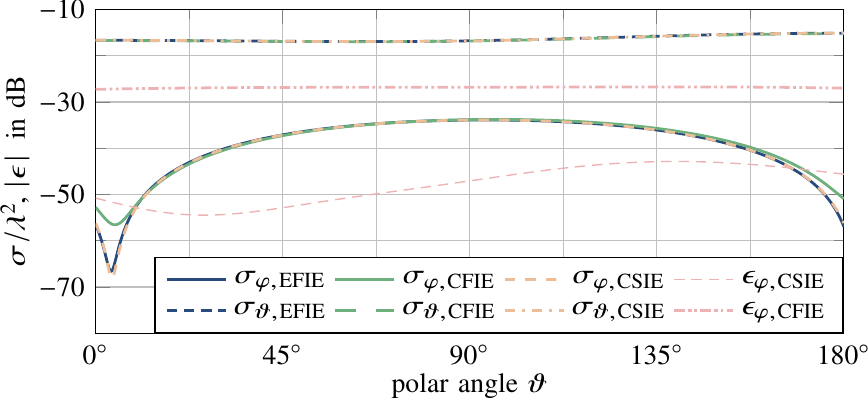}%
  }%
 \end{minipage}
\caption{Bistatic RCS for the small wedge,  with different weighting factors in the CFIE and CSIE. (a) Comparison of $0.1$-CFIE and M-CSIE  in the $\vartheta=90^\circ$ plane. (b) Comparison of CFIE and CSIE  in the $\varphi=90^\circ$ plane.\label{fig:RCS_wedge}}
\end{figure}

The bistatic radar cross section ({RCS}) $\sigma$ has been computed in the far field for an incident plane wave from $-$\textit{z} direction with $E_x$ linear polarization. 
The resulting scattered electric field in the $\vartheta=\uppi/2$ plane is given in Fig.~\ref{fig:RCS_wedge}(a), for the comparison of $0.1$-CFIE and M-CSIE. 
The {$0.1$-CFIE} shows a relative error of up to $-12$\,dB compared to the standard EFIE solution, with the error defined as
\begin{equation}
\epsilon_i = \left|\frac{E_{i,\mathrm{ref}}-E_{i}}{\max\limits_{\vartheta, \varphi} \left(E_{\vartheta}(\vartheta,\varphi),E_{\varphi}(\vartheta,\varphi)\right)}\right|.
\end{equation}
The M-CSIE shows a comparatively lower error of below $-24$\,dB.  
In Fig.~\ref{fig:RCS_wedge}(b), the scattered electric far field is given in the $\varphi=\uppi/2$ plane. 
It can be seen that the CFIE shows an error level  of about $-26$\,dB, whereas the CSIE shows an error of below $-42$\,dB, compared to the EFIE solution, which is in the range of the achievable accuracy for the chosen low-order discretization.

\subsection{Pyramid}

\begin{figure}[tp]
\includegraphics{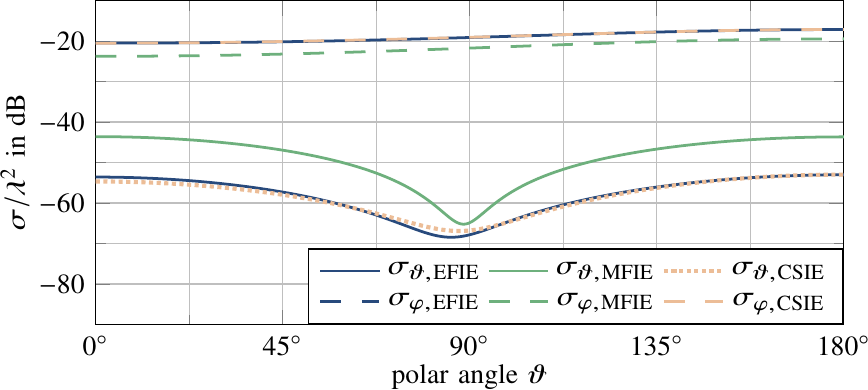}
\caption{Bistatic RCS for the pyramid in the $\varphi=90^\circ$ plane.\label{fig:RCS_pyramid}}
\end{figure}

As a second example, a small pyramid with 36 electric (and magnetic) current unkowns is considered. 
The meshed structure owning 24 triangular cells is depicted in Fig.~\ref{fig:structures}(b). 
The mean edge length is $\lambda/10$, the height in $z$ direction is $0.19\lambda$ and the base square has an edge length of $0.135\lambda$. 
Again, the  EFIE, MFIE, CFIE and CSIE solutions for the scattering problem are compared, with a  plane wave incident from $-$\textit{z} direction and $E_x$ polarization. 
The results for the bistatic RCS are given in Fig.~\ref{fig:RCS_pyramid} in the  $\varphi=\uppi/2$ plane. 
The MFIE shows an error of $-13$\,dB and the CFIE shows an error of $-28$\,dB as compared to the EFIE. 
For the proposed CSIE formulation, the error level is below $-41$\,dB. 
These two examples, the wedge and the pyramid, clearly show that the weak CSIE formulation can circumvent the accuracy problem of the MFIE and also of the CFIE.

\subsection{Flamme}

\begin{figure}[tp]
\centering
\includegraphics{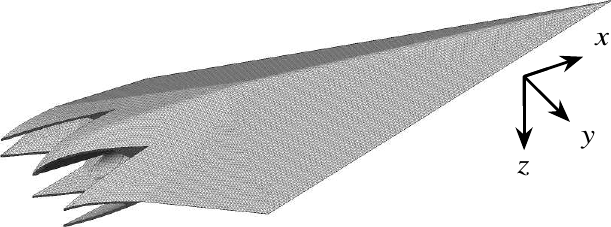}
\caption{Meshed model of the stealth plane Flamme.\label{fig:structure_flamme}}
\end{figure}

A more complicated and realistic example is the stealth object Flamme~\cite{GuerelBagciCastelliEtAl2003,Eibert_2005}. 
The structure is discretized with a $\lambda/10$ mesh, which is depicted in Fig.~\ref{fig:structure_flamme}. 
This means $79173$ unkown electric  currents on 52782 triangular mesh cells, and in addition the same amount of magnetic current unkowns for the CSIE. 
Again, the bistatic RCS is computed for an incident plane wave from an incident angle of $\varphi=80\si{\degree}$. 
\begin{figure}[tp]
\centering
\includegraphics{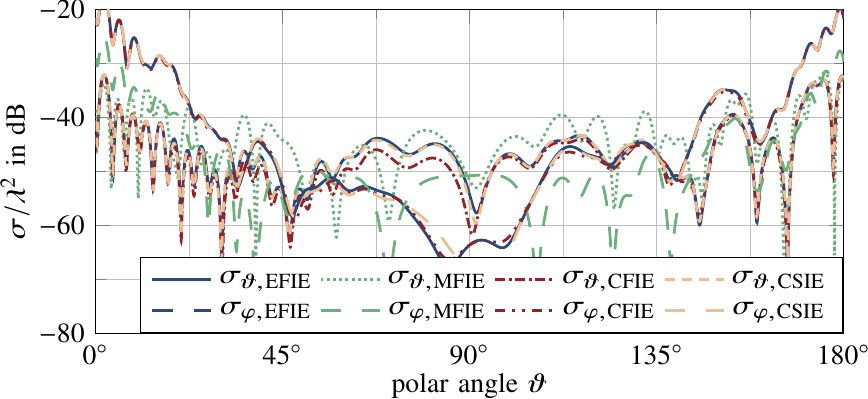}
\caption{Bistatic RCS for Flamme with a $\lambda/10$ mesh at 10\,GHz, in the $\varphi=0^\circ$ plane.\label{fig:flamme_RCS}}
\end{figure}
The MFIE exhibits an error of $-30$\,dB, the CFIE an error of $-40$\,dB and the CSIE an error of $-52$\,dB as compared to the EFIE solution. 

Furthermore, the structure is simulated at a lower frequency, which then results in a $\lambda/25$ mesh. 
The scattered electric far field  is given in Fig.~\ref{fig:flamme_RCS}  in the $\varphi=0^\circ$ plane.
\begin{figure}[tp]
\centering
\includegraphics{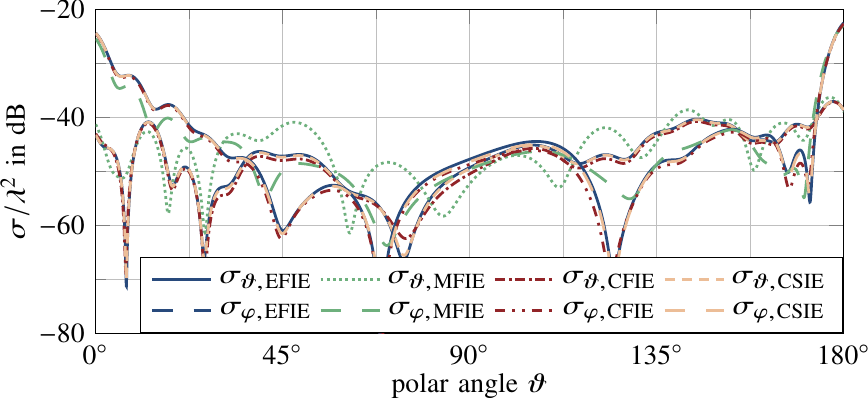}
\caption{Bistatic RCS for Flamme with a $\lambda/25$ mesh at 4\,GHz, in the $\varphi=0^\circ$ plane.\label{fig:flamme_RCS}}
\end{figure}
Compared to the EFIE solution, the MFIE has an accuracy of $-21$\,dB only, possibly due to an interior resonance and also the sharp edges. 
The CFIE circumvents interior resonance problems 
and shows a smaller error of $-34$\,dB. 
The CSIE shows the best accuracy with an error level of  $-44$\,dB as compared to the EFIE solution.

\begin{figure}[tp]
\centering
\includegraphics{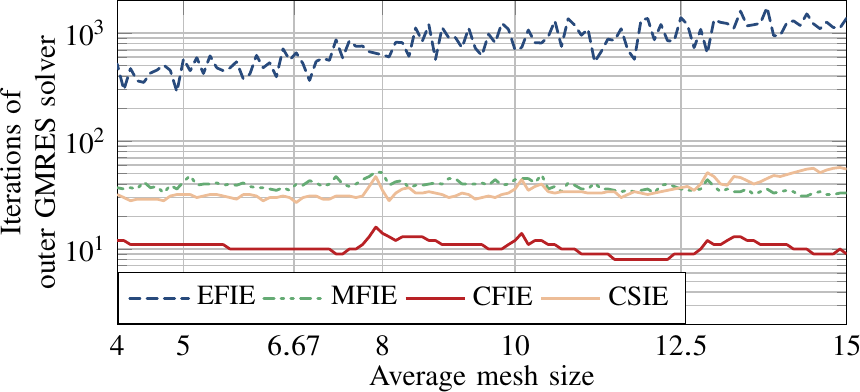}
\caption{Comparison of EFIE, MFIE, CFIE and CSIE iterative solver convergence to a residual error of $10^{-4}$ for the object Flamme.}\label{fig:convergence}
\end{figure}

For the iterative solution of the MoM equation system, a generalized minimal residual (GMRES) solver was employed in an inner-outer scheme. 
The iterative solver convergence for the outer loop is analyzed in Fig.~\ref{fig:convergence} for the same mesh at different frequencies. 
 To obtain a residual error of $10^{-4}$, CSIE and MFIE converge within about the same number of iterations. 
As expected, the CFIE is even faster. 
The worse conditioning of the EFIE is obvious, as it converges much slower than the CSIE. 
This is especially remarkable, as only the CSIE shows the same accuracy as the EFIE. 

\section{Conclusion}
As the adjoint formulation of the combined field integral equation (CFIE), the combined source integral equation (CSIE) is able to resolve the interior resonance problem of the electric field integral equation (EFIE) for perfect electrically conducting (PEC) objects. 
In the presented approach with a weak form combined source condition, the accuracy of the solution is improved as compared to MFIE  and also CFIE solutions of PEC scattering problems.  
The main advantage of the proposed method is the use of only Rao-Wilton-Glisson (RWG) basis and testing functions, where the conversion of the primary unknown electric surface current densities into magnetic surface current densities is achieved by a weak form side condition, namely a modified weak form Leontovich impedance boundary condition (IBC). 
The magnetic current unknowns double the dimensionality, but the memory consumption is only increased negligibly. 
The iterative solver convergence of the presented CSIE formulation is faster than for the EFIE, while obtaining the same accuracy in the simulation results. 
Therefore, the proposed CSIE formulation  is found to be a very promising alternative to the CFIE, with regard to solution accuracy and also implementational and computational effort with RWG functions only.

\ifCLASSOPTIONcaptionsoff
  \newpage
\fi

\bibliographystyle{IEEEtran}
\bibliography{ref2}
\end{document}